\documentclass[preprint]{elsarticle}

\usepackage{graphicx}% Include figure files
\usepackage{dcolumn}% Align table columns on decimal point
\usepackage{bm}% bold math
\usepackage{epsfig}
\usepackage{subfigure}
\usepackage{amssymb}

\begin{document}

\title{Constrained Minimum-Energy Optimal Control of the Dissipative Bloch Equations}
%\title{Minimum-Energy Optimal Control of Bloch Equations in the case of Dominant Transverse Relaxation and with Bounded Control Amplitude}

\author[a]{Dionisis Stefanatos}
\ead{dionisis@post.harvard.edu}

\author[b]{Jr-Shin Li}
\ead{jsli@seas.wustl.edu}

\address[a]{Prefecture of Kefalonia, Argostoli, Kefalonia 28100, Greece}
\address[b]{Washington University, St. Louis, MO 63130, USA}

\begin{abstract}
In this letter, we apply optimal control theory to design minimum-energy $\pi/2$ and $\pi$ pulses for the Bloch system in the presence of relaxation with constrained control amplitude. We consider a commonly encountered case in which the transverse relaxation rate is much larger than the longitudinal one so that the latter can be neglected. Using the Pontryagin's Maximum Principle, we derive optimal feedback laws which are characterized by the number of switches, depending on the control bound and the coordinates of the desired final state.

%In this letter, we apply optimal control theory to design minimum-energy $\pi/2$ and $\pi$ pulses for Bloch equations, in the commonly encountered case where transverse relaxation rate is much larger than longitudinal, so the latter can be neglected, and with constrained control amplitude. Using Pontryagin's Maximum Principle we derive an optimal feedback law which, depending on the control bound and the coordinates of the desired final point, displays zero, one or two switchings.
\end{abstract}


\maketitle

\section{Introduction}

Optimal control theory \cite{Pontryagin} has been extensively used recently for the design of pulses that optimize the performance of various Nuclear Magnetic Resonance (NMR) and quantum systems limited by the presence of relaxation \cite{Khaneja03_1,Khaneja03_2,BBCROP,TROPIC,Stefanatos04,Stefanatos05,Sugny,Bonnard09,BonnardIEEE,Wang, Stefanatos09,Li}, the dissipation due to random interactions between the system and its environment. In this letter, we employ tools from optimal control to derive minimum-energy $\pi/2$ and $\pi$ pulses for a simple NMR system described by the Bloch equations. In particular, we study the case where transverse relaxation dominates the dissipation of the system and the control amplitude is bounded.

The problem of optimal control of the Bloch equations and its closely related corresponding problem for a two-level quantum system have received considerable attention. D' Alessandro and Dahleh \cite{D'Alessandro} considered the problem of minimum-energy optimal control for a two-level quantum system without dissipation. Boscain and Mason \cite{Boscain} examined the time minimal problem for a spin-$1/2$ particle in a magnetic field neglecting dissipation. Sugny, Kontz and Jauslin \cite{Sugny}, Bonnard and Sugny \cite{Bonnard09} and Bonnard, Chyba and Sugny \cite{BonnardIEEE} studied extensively the problem of time-optimal control for a dissipative two-level quantum system.

In our recent work, we studied the problem of designing minimum-energy $\pi/2$ and $\pi$ pulses for the Bloch system dominated by the transverse relaxation with unlimited control amplitude \cite{Stefanatos09}. This dissipative system is of great practical importance as it is a very good approximation for many applications of interest. In this article, we extend this previous work to consider the case where the control amplitude is limited, which accounts for realistic limitations of the experimental setup and also makes the problem more interesting from a control theoretic perspective.

%In the current letter we still focus on the case of dominant transverse relaxation, but here an upper bound is imposed on the control. The control bound accounts for realistic limitations of the experimental setup and also makes the problem more interesting from a control-theoretical perspective, since it is now necessary to use Maximum Principle rather than Calculus of Variations for the solution}.

In the next section, we formulate the related optimal control problems of such pulse designs. The solutions of these problems are presented in Section 3, which is the main contribution of this article. Then in Section 4, we present some examples to demonstrate our analytical results.

%In the next section we formulate the optimal control problem, which is solved in section 3, the core of this letter. In section 4 we present some examples while section 5 concludes the article.

%%%%%%%%%%%%%%%%%%%%%%%%%%%%%%%%%%%%%%%%%%%%%%%%%%%%%%%%%%%
\section{Optimal Control of Dissipative Bloch Systems}

\begin{figure}[t]
 \centering
		\begin{tabular}{cc}
     	    \subfigure[$\ $$u \perp \mathbf{r}_\perp$]{
	            \label{fig:3d}
	            \includegraphics[width=.45\linewidth]{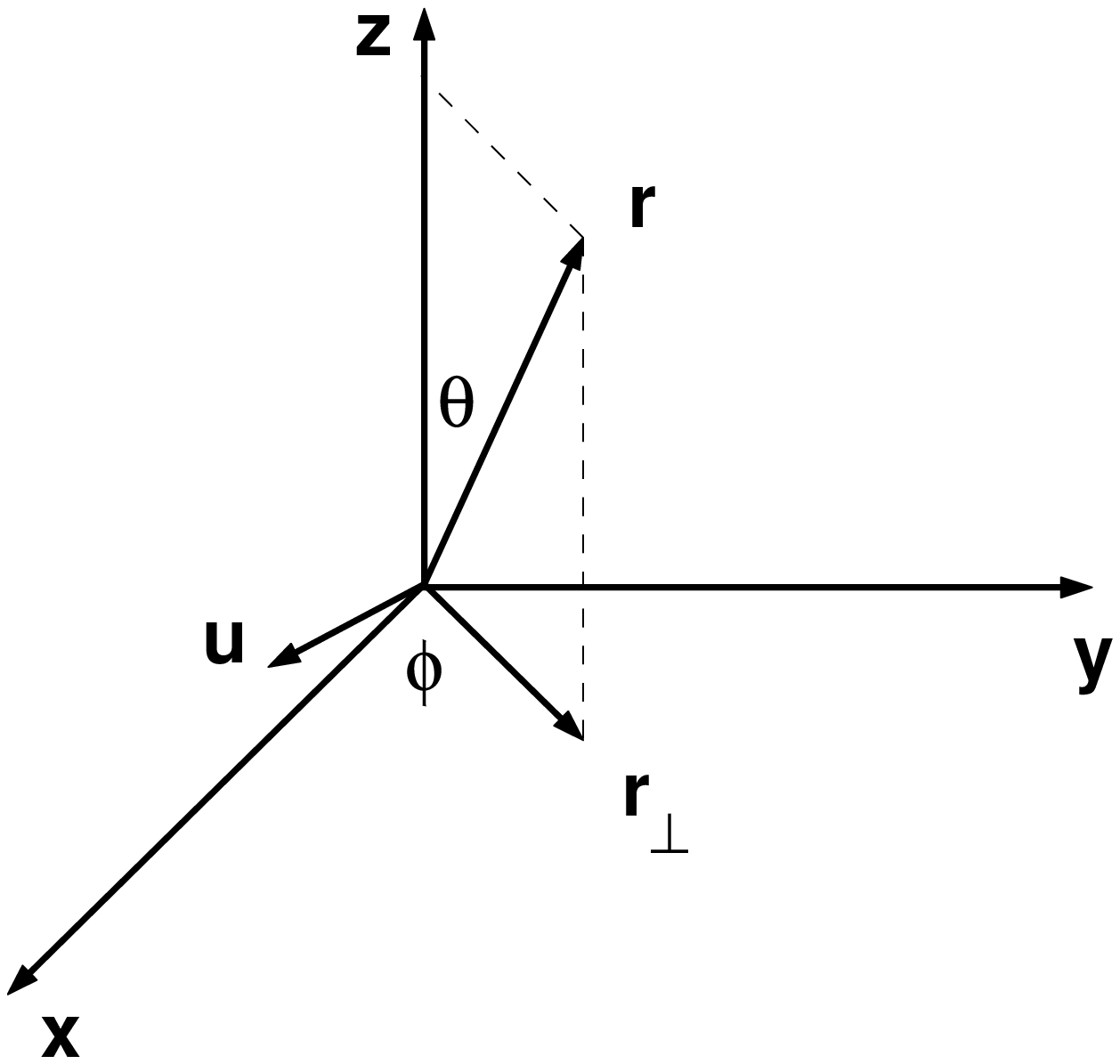}} &
	        \subfigure[$\ $$u \parallel x$-axis]{
	            \label{fig:2d}
	            \includegraphics[width=.45\linewidth]{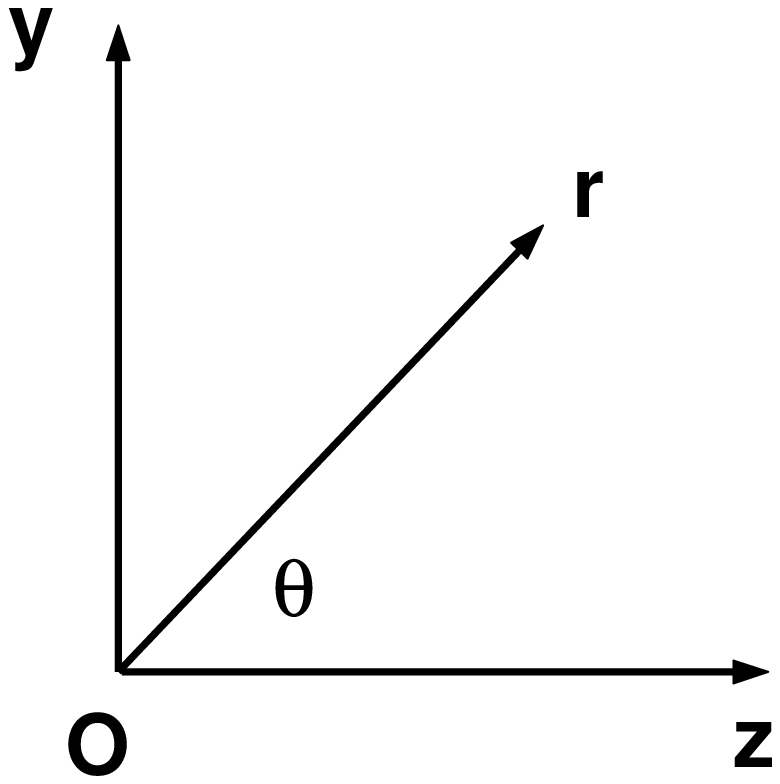}}
		\end{tabular}
\caption{The optimal transverse magnetic field $u$ is perpendicular to $\mathbf{r}_\perp$ and its phase is constant (panel a). Without loss of generality, the experimental setup can be arranged such that $u \parallel x$-axis. In this case, $\phi=\pi/2$ and $\mathbf{r}$ rotates in $yz$-plane (panel b). For convenience we display $z$ in the horizontal axis and $y$ in the vertical, so that $(r,\theta)$ have the common configuration of polar coordinates on the plane.}
\label{fig:spherical}
\end{figure}

%\subsection{Problem Definition}
In a resonant rotating frame, the Bloch equations with the longitudinal relaxation neglected are of the form \cite{Ernst}
\begin{eqnarray}
	\dot{z} & = & u_yx-u_xy \\
	\dot{x} & = & -Rx-u_yz \\
    \label{y}\dot{y} & = & -Ry+u_xz,
\end{eqnarray}
where $\mathbf{r}=(x, y, z)$ is the magnetization vector, $u_x,u_y$ are the transverse components of the magnetic field and $R>0$ is the transverse relaxation rate. The above equations constitute a dissipative bilinear control system. By the following change of variables (see Fig. \ref{fig:3d}) and time rescaling
\begin{eqnarray}
	a & = & \ln r=\ln(\sqrt{x^2+y^2+z^2}) \nonumber\\
	\tan\theta & = & \sqrt{x^2+y^2}/z \nonumber\\
    \tan\phi & = & y/x \nonumber\\
    t_{new} & = & R\,t_{old},\nonumber
\end{eqnarray}
we arrive at a new system
\begin{eqnarray}
	\label{magnitude}\dot{a} & = & -\sin^2\theta \\
	\label{angle}\dot{\theta} & = & u-\sin\theta\cos\theta \\
    \label{phi}\dot{\phi} & = & v\cot\theta
\end{eqnarray}
where $u = (u_x/R)\sin\phi-(u_y/R)\cos\phi, v = (u_x/R)\cos\phi+(u_y/R)\sin\phi$ are the normalized components of transverse magnetic field perpendicular and parallel to $\mathbf{r}_\perp = (x,y)$, respectively.

Note that $v$ does not affect the angle $\theta$ of the pulse. It just rotates $\mathbf{r}$ around $z$-axis, resulting in a waste of energy. Thus, optimality requires that $v=0$ and hence $\phi=\mbox{constant}$. If for example we choose $u$ to be parallel to the $x$-axis, then $\phi=\pi/2$ and $\mathbf{r}$ rotates in $yz$-plane, see Fig. \ref{fig:2d}. This is the case that we consider in this letter. Equations (\ref{magnitude}) and (\ref{angle}) are sufficient to describe this rotation.

%% =================== Figure 1 ==================

%Observe that when $\omega$ is unbounded, we can rotate $\theta$ instantaneously to the desired final value $\theta(T)=\pi/2\,\mbox{or}\,\pi$ without losses in $a$, i.e. with $a(T)=a(0)=0$ (equivalently $M(T)=M(0)$). This transfer requires an infinite amount of energy so it is unrealistic.

The optimal control problem that we would like to pursue is formulated as follows. Consider the dynamical system as in (\ref{magnitude}), (\ref{angle}), starting from $(r(0),\theta(0))=(1,0)$ (corresponding to $(a(0),\theta(0))=(0,0)$) and for a specified final value $a(\tau)<a(0)=0$ (equivalent to $r(\tau)<r(0)=1$), what is the optimal bounded control $u(t)$ with $|u(t)|\leq m\in\mathbb{R}^+$, $0\leq t\leq \tau$, that accomplishes the transfer between the above starting point and the target point $(a(\tau),\theta(\tau)=\pi/2\,\mbox{or}\,\pi)$, while minimizing the energy $E=\int_0^\tau u^2(t)/2dt$? Note that for the transfers that we study here, it must be $m>1/2$, otherwise Eq. (\ref{angle}) reaches an equilibrium point $\theta_0<\pi/2$. Also, the final time $\tau$ is unspecified.

%The starting point is $(r(0),\theta(0))=(1,0)$, corresponding to $(a(0),\theta(0))=(0,0)$. In this letter we solve the following optimization problem: for a specified final value $a(\tau)<a(0)=0$ (equivalently specified $r(\tau)<r(0)=1$), what is the optimal bounded control $u(t)$, $|u(t)|\leq m$, $0\leq t\leq \tau$, that accomplishes the transfer from the above starting point to the final point $(a(\tau),\theta(\tau)=\pi/2\,\mbox{or}\,\pi)$, while minimizing energy $E=\int_0^\tau u^2(t)/2dt$? Note that for the transfers that we study here, it must be $m>1/2$, otherwise Eq. (\ref{angle}) reaches an equilibrium point $\theta_0<\pi/2$. Also, the final time $\tau$ is unspecified.

%%%%%%%%%%%%%%%%%%%%%%%%%%%%%%%%%%%%%%%%%%%%%%%%%%%%%%%%%
\section{Derivation of the Optimal Control}
The control Hamiltonian for the addressed problem is
\begin{equation}
\label{hamiltonian}
H = -u^2/2+\lambda_{\theta}(u-\sin\theta\cos\theta)-\lambda_a\sin^2\theta
\end{equation}
where $\lambda_{\theta},\lambda_a$ are the Lagrange multipliers. According to Pontryagin's maximum principle \cite{Pontryagin}, the necessary conditions for optimality of $(u(t),\theta(t), a(t), \lambda_{\theta}(t),\lambda_a(t))$ are
\begin{eqnarray}
\label{adjoint1}
\dot{\lambda}_\theta & = & -\partial H/\partial\theta = \lambda_\theta \cos 2\theta+\lambda_a \sin 2\theta\\
\label{adjoint2}
\dot{\lambda}_a & = & -\partial H/\partial a = 0 \\% \Rightarrow \lambda_a=\mbox{constant}\\
\label{optimality1}
u & = & \mbox{arg}\max_u H(u,\theta,a,\lambda_\theta,\lambda_a)
\end{eqnarray}
From (\ref{adjoint2}), we immediately know that $\lambda_a$ is a constant. Additionally, the optimal $(u,\theta,a, \lambda_{\theta},\lambda_a)$ satisfies \cite{Pontryagin}
\begin{equation}
\label{hzero}
H(u,\theta,a,\lambda_\theta,\lambda_a)=0,\; 0\leq t\leq \tau.
\end{equation}

Let $E=\min_u\int_0^\tau u^2(t)/2dt$ be the minimum cost corresponding to the optimal solution. Using calculus of variations we find that small changes in the $a$-coordinate of the final point $\delta a_\tau$, and small changes in the final time $\delta \tau$, produce the following change in the minimum cost
\begin{equation}
\label{variation}
\delta E=\lambda_a(\tau)\,\delta a_\tau-H(\tau)\,\delta \tau.
\end{equation}
Therefore
\begin{equation}
\label{dynamic}
\lambda_a(\tau)=\partial E/\partial a_\tau=\partial E/\partial r_\tau \cdot dr_\tau/da_\tau=\partial E/\partial r_\tau\cdot r_\tau,
\end{equation}
where $r_\tau=e^{a_\tau}$ is the radius of the final point.

%% =================== Figure 2 ==================
\begin{figure}[t]
 \centering
		\begin{tabular}{cc}
     	\subfigure[$\ $Case of $\pi$ pulse]{
	            \label{fig:sign1}
	            \includegraphics[width=.45\linewidth]{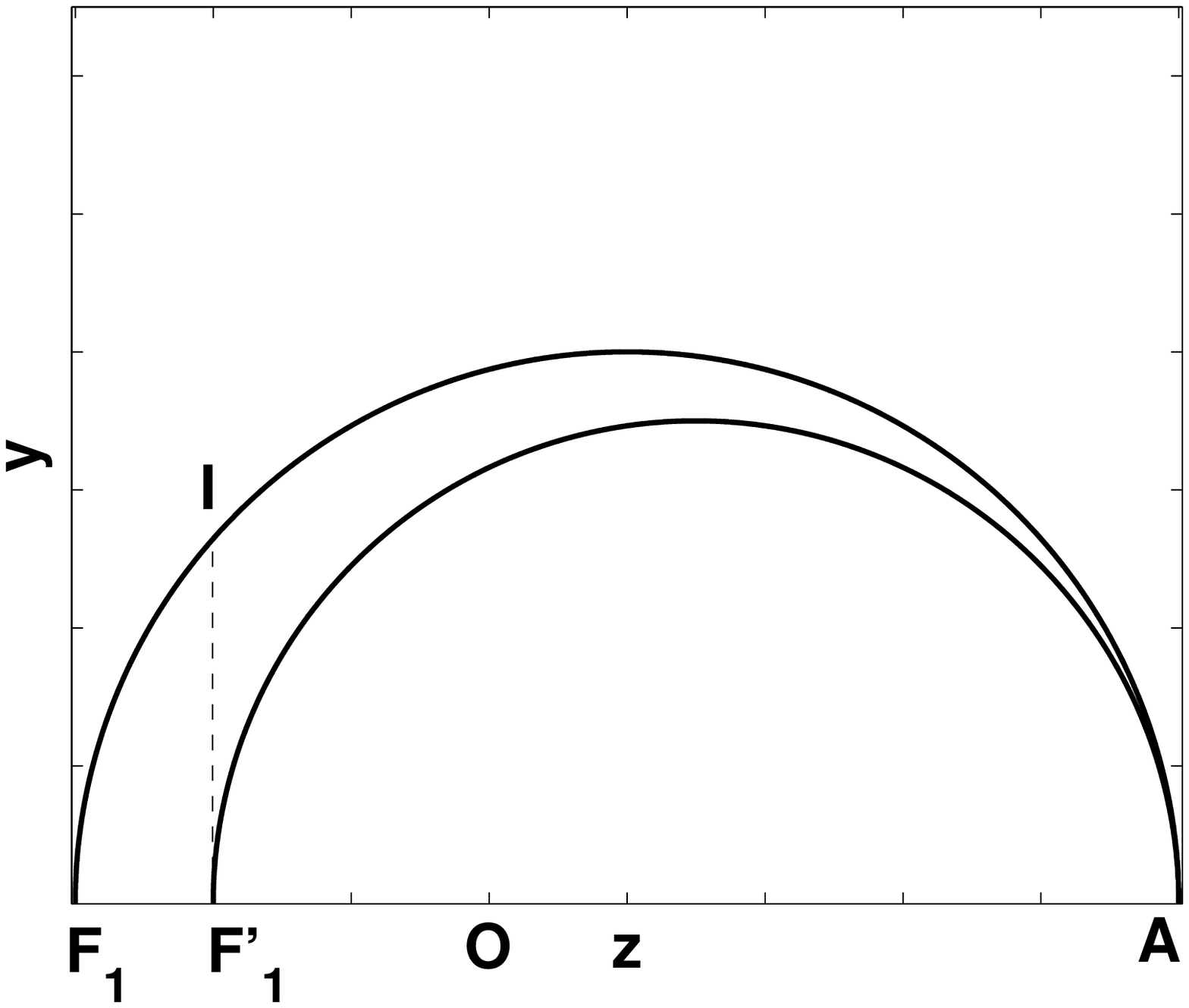}} &
	        \subfigure[$\ $Case of $\pi/2$ pulse]{
	            \label{fig:sign2}
	            \includegraphics[width=.45\linewidth]{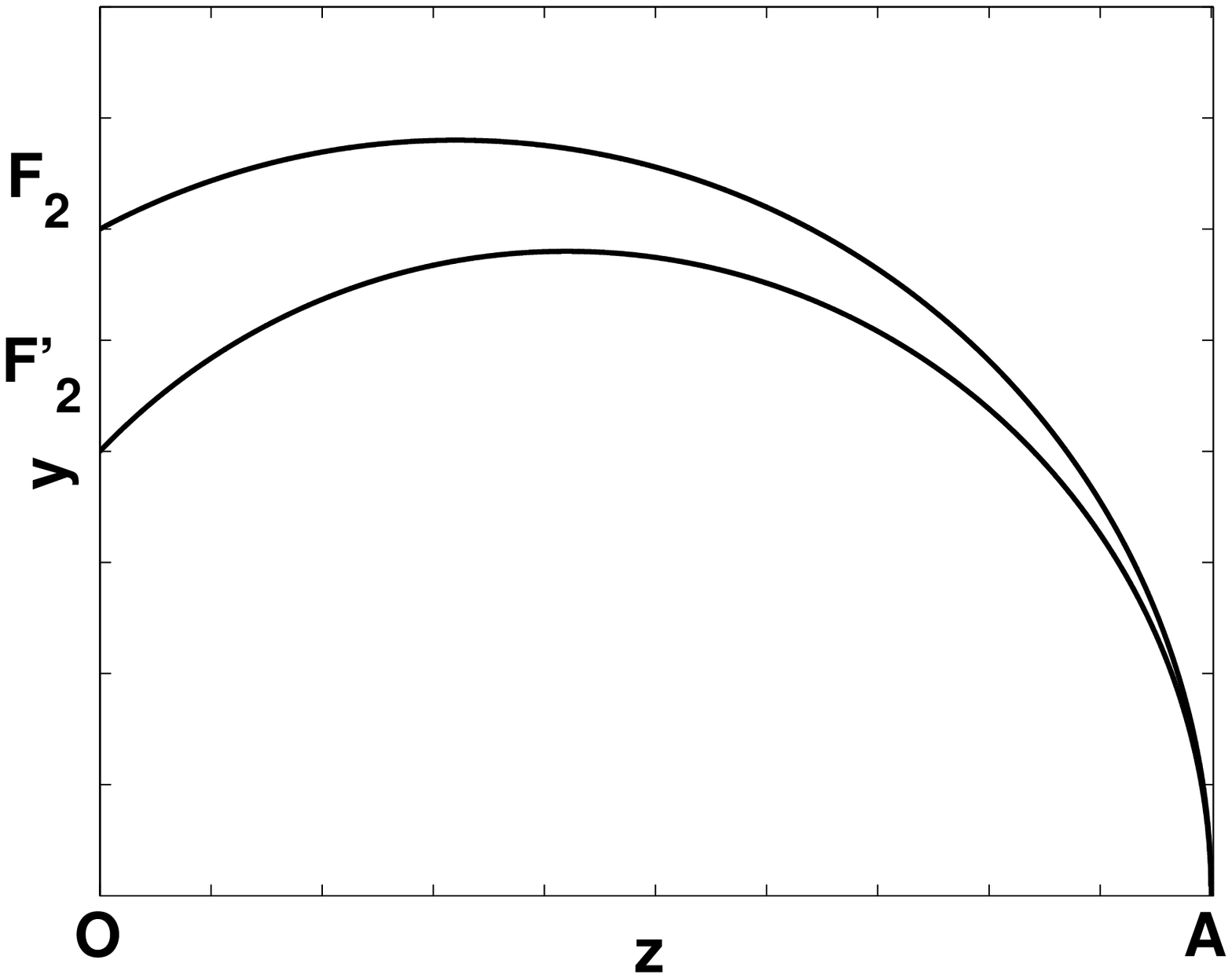}}
		\end{tabular}
\caption{The final point $F'_1$ can be reached by either following the minimum-energy path $AF'_1$, or traveling along $AF_1$ up to point $I$ and then leaving the system relax to $F'_1$ (panel a). Analogously, $F'_2$ can be reached by either following the minimum-energy path $AF'_2$, or traveling along $AF_2$ and then leaving the system relax to $F'_2$ (panel b). .}
\label{fig:Evsr}
\end{figure}

For the transfers that we examine here it is
\begin{equation}
\label{sign}
\partial E/\partial r_\tau\geq 0
\end{equation}
i.e. the larger the final $r$, the more energy is needed to quickly rotate the vector before it dissipates. To see this, we refer to Fig. \ref{fig:Evsr}. Let $E, E'$ be the minimum energies necessary to reach the final points $F_1(r_\tau,\pi),F'_1(r'_\tau,\pi)$, respectively, with $r'_\tau\leq r_\tau$, see Fig. \ref{fig:sign1}. The corresponding minimum-energy paths are $AF_1,AF'_1$. An alternative way to reach the point $F'_1$ is the following: travel along $AF_1$ up to the point $I$, with $z_I=z_{F'_1}=-r'_\tau$, then set $u=0$ and wait until dissipation eliminates the $y$-coordinate, see (\ref{y}). The energy $E''$ spent for this travel is the portion of $E$ necessary to reach $I$, so $E''\leq E$. By the definition of $E'$ it is also $E'\leq E''$, and thus $E'\leq E$. Analogously, let $E, E'$ be the minimum energies necessary to reach the final points $F_2(r_\tau,\pi/2),F'_2(r'_\tau,\pi/2)$, respectively, with $r'_\tau\leq r_\tau$ again, see Fig. \ref{fig:sign2}. The corresponding minimum-energy paths are $AF_2,AF'_2$. An alternative way to reach the point $F'_2$ is the following: travel along $AF_2$ up to the point $F_2$, then set $u=0$ and wait until dissipation brings the system at the point $F'_2$, see (\ref{y}). The energy spent for this travel is the necessary energy to reach $F_2$, i.e. $E$. Then, by definition of $E'$, it is $E'\leq E$. Thus (\ref{sign}) is true and from (\ref{dynamic}) we have that $\lambda_a(\tau)\geq 0$. But $\lambda_a=\mbox{constant}$, so we can set
\begin{equation}
\label{lambdaa}
\lambda_a=\kappa^2/2,\,\,\kappa\geq 0.
\end{equation}

Having determined the sign of $\lambda_a$, we first examine the case of unbounded control and then we use the developed intuition to study the general case of bounded control.

\subsection{Unbounded Control}
When the control $u$ is unbounded, then from (\ref{optimality1}) we conclude that $\partial H/\partial u = 0$. This condition and Eq. (\ref{hamiltonian}) yield
\begin{equation}
\label{optimality}
u = \lambda_{\theta}.
\end{equation}
%\begin{equation}
%\label{optimality}
%\partial H/\partial u = 0 \Rightarrow u = \lambda_{\theta}.
%\end{equation}
Using (\ref{optimality}) and (\ref{lambdaa}), the condition (\ref{hzero}) becomes
\begin{equation}
\label{optimallambda}
\lambda^2_\theta-2\lambda_\theta \sin\theta\cos\theta-\kappa^2 \sin^2\theta=0.
\end{equation}
The optimal $u$ is then given by the following feedback law
\begin{equation}
\label{u}
u(\theta)=\lambda_{\theta}=\sin\theta(\cos\theta+\sqrt{\cos^2\theta+\kappa^2}).
\end{equation}
Note that only the positive solution of the quadratic equation has physical meaning (corresponds to increasing $\theta$) for the transfers that we study here. Using (\ref{optimality}) and (\ref{u}), the validity of (\ref{adjoint1}) can be easily verified.
Inserting (\ref{u}) in (\ref{angle}) we obtain the differential equation for the optimal trajectory
\begin{equation}
\label{theta}
\dot{\theta}=\sin\theta\sqrt{\cos^2\theta+\kappa^2}
\end{equation}
Eliminating time between (\ref{magnitude}) and (\ref{theta}) we obtain
\begin{equation}
\label{difeq}
\frac{da}{d\theta}=-\frac{\sin\theta}{\sqrt{\cos^2\theta+\kappa^2}}.
\end{equation}
Integrating the above equation from the starting point $(0,0)$ to the point \textbf{$(\ln r,\theta)$}, we find the optimal trajectory
\begin{equation}
\label{opttraj}
r(\theta)=\frac{\cos\theta+\sqrt{\cos^2\theta+\kappa^2}}{1+\sqrt{1+\kappa^2}}
\end{equation}
Setting $\theta_\tau=\pi/2,\pi$ in the above equation, we find the optimal $\kappa$ for the $\pi/2$ and $\pi$ pulses, as a function of the radius $r_\tau$ of the final point
\begin{equation}
\label{kappa}
\kappa_{\pi/2}=\frac{2r_\tau}{1-r^2_\tau},\;\;\kappa_{\pi}=\frac{2\sqrt{r_\tau}}{1-r_\tau}
\end{equation}
%Note that since $\theta=0,\pi$ are equilibrium points for (\ref{theta}), we actually start from a small positive initial value $\theta(0)=\epsilon$ for both cases,additionally for the $\pi$ pulse we end to the value $\theta(T)=\pi-\epsilon$. The duration of the optimal pulses as $\epsilon\rightarrow 0$

The energy of the optimal pulses is calculated from
\begin{equation}
E=\int_0^\tau\frac{u^2(t)}{2}dt=\int_0^{\theta_\tau}\frac{u^2(\theta)}{2\dot{\theta}(\theta)}\,d\theta
\end{equation}
using (\ref{u}) and (\ref{theta}). The results for $\theta_\tau=\pi/2,\pi$ are
\begin{equation}
\label{energy}
E_{\pi/2}=\frac{1}{1-r^2_\tau},\;\;E_{\pi}=\frac{1+r_\tau}{1-r_\tau}
\end{equation}
Using (\ref{energy}) in (\ref{dynamic}), it is easy to verify the validity of (\ref{lambdaa}) with $\kappa$ given by (\ref{kappa}).
%Observe that for $r\rightarrow 1$ we have $T_{\pi/2}, T_{\pi}\rightarrow 0$ and $E_{\pi/2}, E_{\pi}\rightarrow \infty$, as mentioned above.
%% =================== Figure 3 ==================

\begin{figure}[t]
 \centering
		\begin{tabular}{cc}
     	\subfigure[$\ $Optimal $\pi/2$ pulse]{
	            \label{fig:pi2}
	            \includegraphics[width=.45\linewidth]{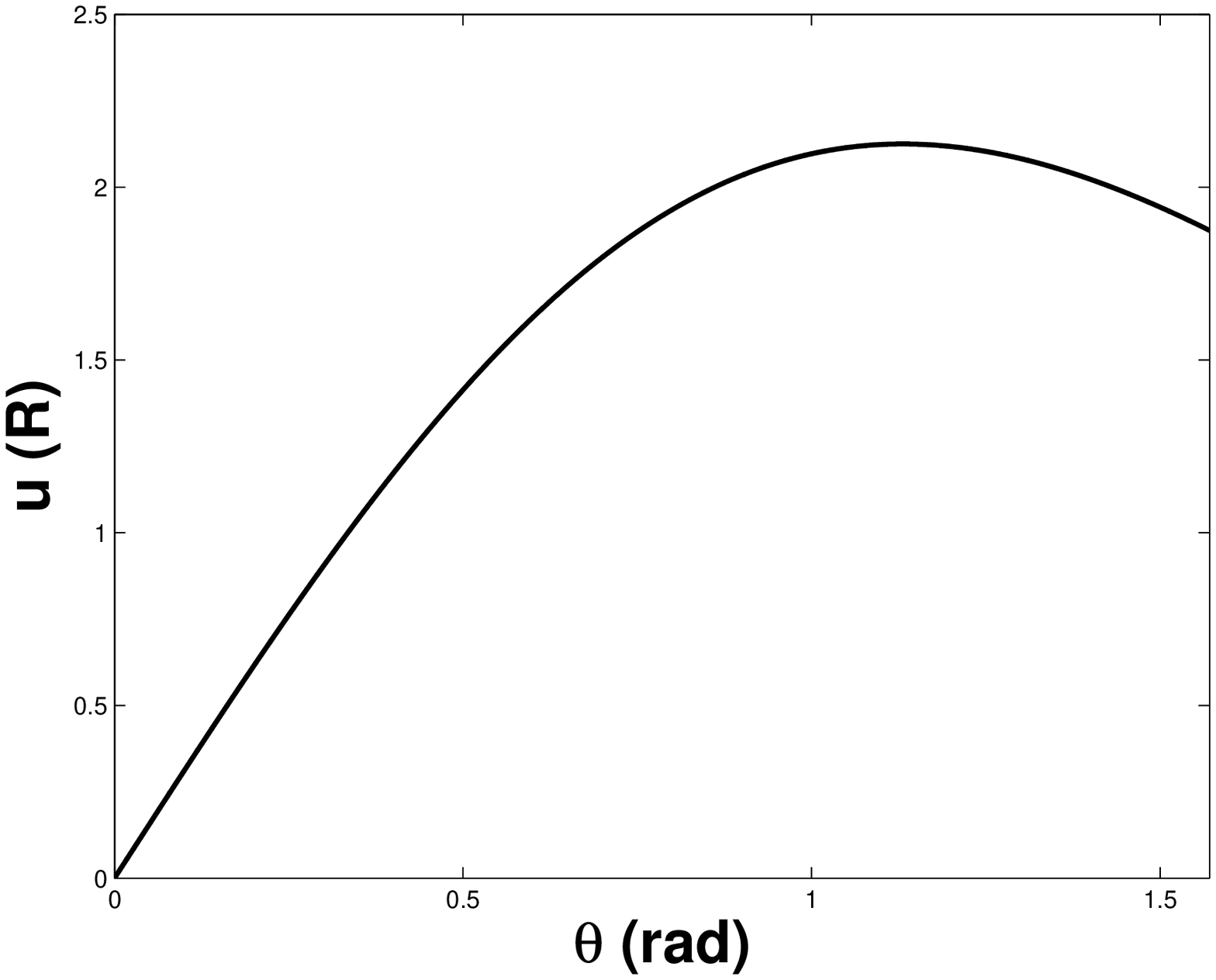}} &
	        \subfigure[$\ $Optimal trajectory]{
	            \label{fig:pi}
	            \includegraphics[width=.45\linewidth]{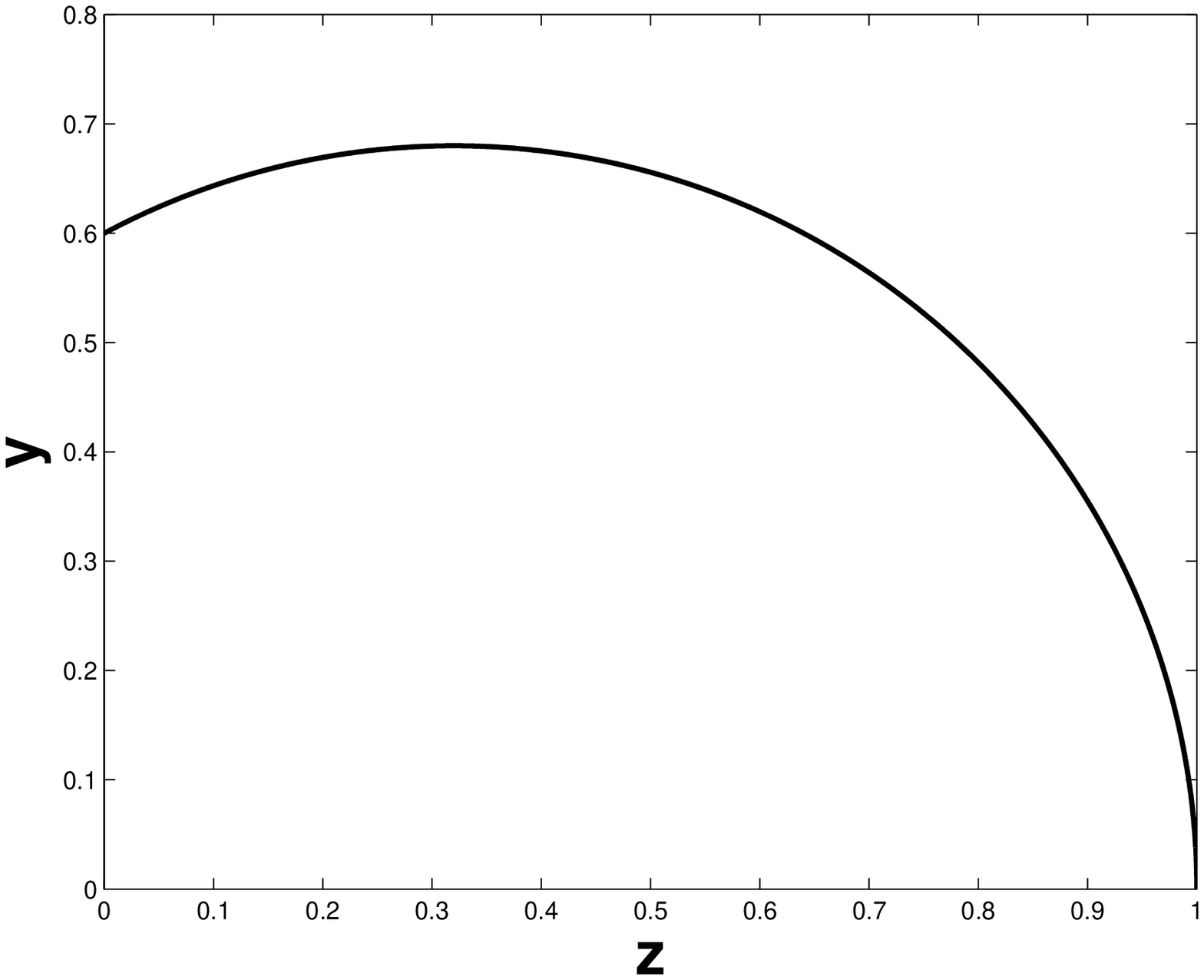}} \\
	        \subfigure[$\ $Optimal $\pi$ pulse]{
	            \label{fig:pi2_traj}
	            \includegraphics[width=.45\linewidth]{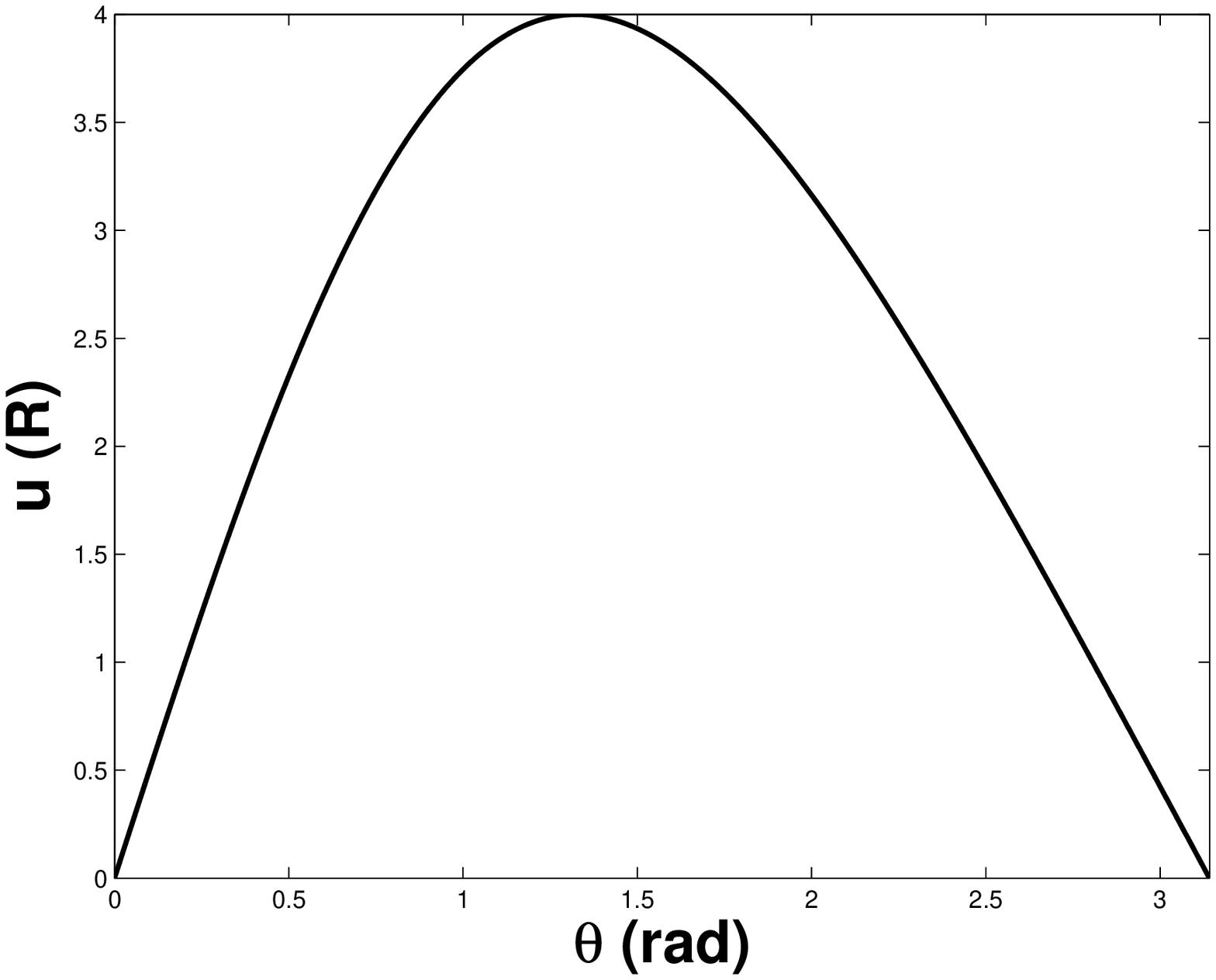}} &
			\subfigure[$\ $Optimal trajectory]{
	            \label{fig:pi_traj}
	            \includegraphics[width=.45\linewidth]{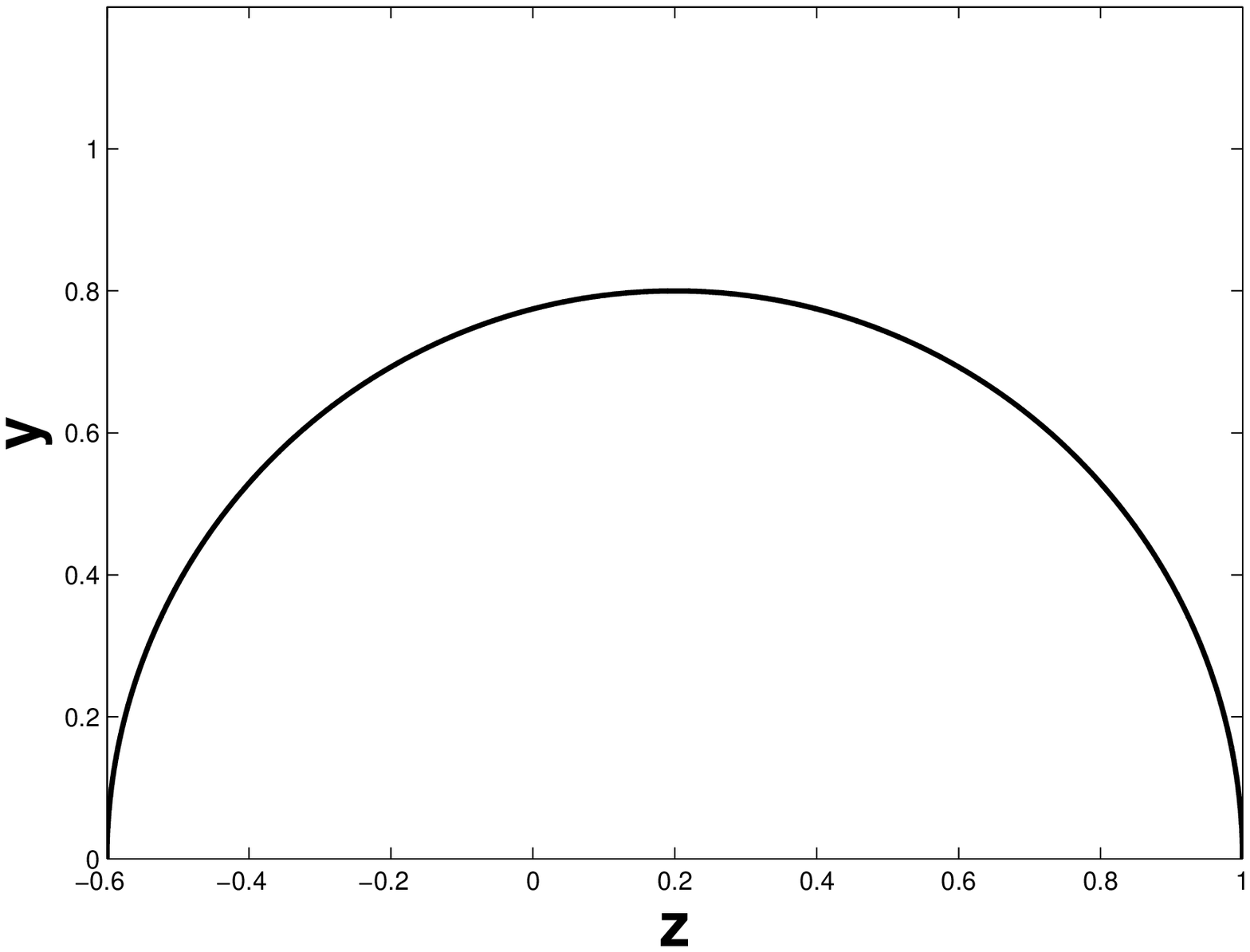}} \\
		\end{tabular}
\caption{minimum-energy $\pi/2$ (panel a) and $\pi$ (panel c) pulses for $r_\tau=0.6$. The corresponding trajectories are also shown (panels b,d).}
 \label{fig:pulses}
\end{figure}
In Fig. \ref{fig:pulses} we plot the optimal $\pi/2$ and $\pi$ pulses for $r_\tau=0.6$, as well as the corresponding trajectories. Observe that optimal $u(\theta)$ is small close to $\theta=0,\pi$ ($z$-axis), directions that are protected against relaxation, while it is large close to $\theta=\pi/2$ ($y$-axis), where dissipation is maximized and thus $\mathbf{r}$ must be rotated faster.

\subsection{Bounded Control}
We now move on to the case where the control is bounded, i.e. $|u|\leq m$, with $m>1/2$ as pointed out before. The control Hamiltonian (\ref{hamiltonian}) is a quadratic form with respect to $u$ that takes its maximum value at $u=\lambda_\theta$, if $|\lambda_\theta|\leq m$, and at the boundary point $u=m$ if $\lambda_\theta>m$. The other boundary point, $u=-m$, corresponds to decreasing $\theta$ and has no physical meaning for the transfers that we examine here.
Initially, the situation is as in the previous case where the optimality condition $u=\lambda_\theta$ holds and the optimal control is given by (\ref{u}). Angle $\theta$ increases following (\ref{theta}) and $u,\lambda_\theta$ change accordingly.
%Now, if at some point the control reaches the maximum allowable value $m$, then a switching takes place and the optimal control takes the value $u(\theta)=m$. The control Hamiltonian (\ref{hamiltonian}), a quadratic form with respect to $u$, takes its maximum value at the border of the allowable interval $|u|\leq m$. This implies that $\lambda_\theta\geq m$ and, as long as this condition holds, the optimal value remains $u(\theta)=m$. In the course of time, $\lambda_\theta$ changes and if crosses the value $m$ to lower values, another switching takes place. During the third phase, the optimal control is given again by (\ref{u}), with the same $\kappa$ as in the first phase, since $\lambda_a=\kappa^2/2$ is constant along the optimal trajectory.
Now suppose that at some point the control reaches the maximum allowable value $m$. From (\ref{u}) we see that this happens at the angles that satisfy the equation
\begin{equation}
\sin\theta(\cos\theta+\sqrt{\cos^2\theta+\kappa^2})=m,
\end{equation}
which is equivalent to the following quadratic equation for $\cot\theta$
\begin{equation}
\label{switch}
\cot^2\theta-\frac{2}{m}\cot\theta+1-\frac{\kappa^2}{m^2}=0.
\end{equation}
If $\theta_1,\theta_2$ are the solutions of (\ref{switch}) in $[0,\pi]$, then it is easy to show that for $\theta\in(\theta_1,\theta_2)$ the following inequality holds
\begin{equation}
\label{ineq}
\sin\theta(\cos\theta+\sqrt{\cos^2\theta+\kappa^2})> m.
\end{equation}
In the interval $\theta\in(\theta_1,\theta_2)$ where the above inequality is true, the relation $u=\lambda_\theta$ gives $u(\theta)>m$ which is not permissible. Thus, the optimal control in this interval is $u(\theta)=m$. From (\ref{hzero}) we can find the Lagrange multiplier $\lambda_\theta$ for the same interval, which is
\begin{equation}
\label{lambda2}
\lambda_\theta(\theta)=\frac{m^2+\kappa^2\sin^2\theta}{2(m-\sin\theta\cos\theta)}
\end{equation}
and $\lambda_\theta(\theta_1)=m$ from (\ref{optimality}).
This $\lambda_\theta$ satisfies the optimality condition (\ref{adjoint1}) with $\theta$ evolving in time according to (\ref{angle}) and $u(\theta)=m$. It is not hard to verify using (\ref{ineq}), (\ref{lambda2}) that $\lambda_\theta> m$ for $\theta\in(\theta_1,\theta_2)$, thus $u(\theta)=m$ is indeed the optimal control in this interval. We call the solutions $\theta_1,\theta_2$ of (\ref{switch}), where a change in the optimal control occurs, the \textit{switching angles}. After the second switching, the optimal control is given again by (\ref{u}), with the same $\kappa$ as in the initial phase, since $\lambda_a=\kappa^2/2$ is constant along the optimal trajectory.

As explained above, the switching angles determine the optimal feedback law.
For $m\geq1$, Eq. (\ref{switch}) has real solutions for $\kappa\geq \sqrt{m^2-1}$.
For $1/2<m<1$ it has real solutions for every $\kappa\geq 0$.
We examine separately these two cases

\subsubsection{$m\geq1$}

%% =================== Figure 4 ==================
\begin{figure}[t]
\centering
\includegraphics[scale=0.45]{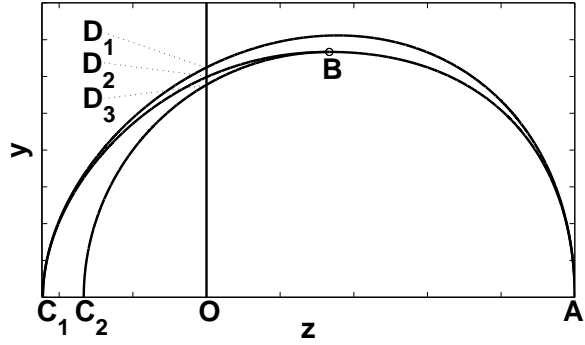}
\caption{Switching curves $AB$ and $BC_1$ for the case $m\geq1$. The outermost curve $AC_1$ defines the reachable set from the starting point $A(1,0)$ when $u\leq m$. Curve $BC_2$ is the optimal trajectory for $\kappa=\sqrt{m^2-1}$ and $\theta\in[\theta_B,\pi]$. The perpendicular axis $\theta=\pi/2$ crosses the curves $AC_1, BC_1, BC_2$ at the points $D_1, D_2, D_3$, respectively.}
\label{fig:mgreater}
\end{figure}

For $\kappa>\sqrt{m^2-1}$ there are two switching angles, given by
\begin{equation}
\label{switchingangle}
\theta_{1,2}=\cot^{-1}\left(\frac{1\pm\sqrt{\kappa^2-m^2+1}}{m}\right)
\end{equation}
In this letter the range of the function $\cot^{-1}$ is considered to be $[0,\pi]$, so
\begin{equation}
\cot^{-1}(x)=\pi-\cot^{-1}(-x),\,\,x<0.
\end{equation}
In the case $\kappa=\sqrt{m^2-1}$, the two angles obtain the common value $\theta_B=\cot^{-1}(1/m)$. The angle of the first switching point is in the range $\theta_{1}\in[0,\theta_B]$, while that of the second $\theta_{2}\in[\theta_B,\pi]$. We find the equation of the first switching curve, i.e. the curve composed by the points $(r_1,\theta_1)$. Before the switching the optimal control is given by (\ref{u}), so each of the points $(r_1,\theta_1)$ belongs to an optimal curve of the form (\ref{opttraj}), i.e.
\begin{equation}
\label{firstswitchingradius}
r_1(\theta_1)=\frac{\cos\theta_1+\sqrt{\cos^2\theta_1+\kappa^2}}{1+\sqrt{1+\kappa^2}}.
\end{equation}
But this $\kappa$ is related to the switching angle $\theta_1$ through (\ref{switchingangle}), which we re-write as
\begin{equation}
\label{kappa1}
\kappa^2=(m\cot\theta_1-1)^2+m^2-1
\end{equation}
since only $\kappa^2$ appears in (\ref{firstswitchingradius}).
These two equations determine the first switching curve, with the angle $\theta_1$ in the interval $\theta_{1}\in[0,\theta_B]$.

After the first switching, the optimal control takes the value $u(\theta)=m$. It maintains this value until the second switching, at angle $\theta_2$. Between the two switchings the evolution is given by
\begin{equation}
\frac{da}{d\theta}=-\frac{\sin^2\theta}{m-\sin\theta\cos\theta},
\end{equation}
as we derive from (\ref{magnitude}), (\ref{angle}) with $u=m$.
Integrating the above equation from $\theta_1$ to $\theta_2$ we find the radius of the second switching point
\begin{equation}
\label{secswitchingradius}
r_2(\theta_2)=r_1(\theta_1)\sqrt{\frac{2m-\sin{2\theta_1}}{2m-\sin{2\theta_2}}}\exp{\left[-\frac{f(\theta_1,\theta_2)}{\sqrt{4m^2-1}}\right]}
\end{equation}
where
\begin{equation}
\label{exponent}
f(\theta_1,\theta_2)=\cot^{-1}\left(\frac{2m\cot\theta_2-1}{\sqrt{4m^2-1}}\right)-\cot^{-1}\left(\frac{2m\cot\theta_1-1}{\sqrt{4m^2-1}}\right)
\end{equation}
and
\begin{equation}
\label{secswitchingangle}
\theta_2=\cot^{-1}(2/m-\cot\theta_1).
\end{equation}
Therefore, to every first switching point $(r_1,\theta_1)$ corresponds a second switching point $(r_2,\theta_2)$ with angle given by (\ref{secswitchingangle}) and radius given by (\ref{secswitchingradius}), (\ref{exponent}). These points compose the second switching curve. The two switching curves are plotted in Fig. \ref{fig:mgreater}, curves $AB$ and $BC_1$. The joint point is $B(r_B,\theta_B)$, where
\begin{equation}
r_B=\frac{\sqrt{m^2+1}}{m+1}.
\end{equation}

The outermost curve $AC_1$ corresponds to the trajectory traveled for $u(\theta)=m, \theta\in[0,\pi]$. The equation for this trajectory can be found by setting $(r_1,\theta_1)=(1,0)$ (starting point $A$) in (\ref{secswitchingradius}), (\ref{exponent}). It is
\begin{equation}
\label{reachable}
r_3(\theta_3)=\sqrt{\frac{2m}{2m-\sin{2\theta_3}}}\exp{\left[-\frac{1}{\sqrt{4m^2-1}}\cot^{-1}\left(\frac{2m\cot\theta_3-1}{\sqrt{4m^2-1}}\right)\right]}
\end{equation}
where we used $(r_3,\theta_3)$ to denote a point on the curve and $\theta_3\in[0,\pi]$. The points between this curve and the horizontal axis define the reachable set from $A$ for a specific control bound $m$. This curve meets the axes $\theta=\pi,\pi/2$ at the points
$C_1(r_{C_1},\pi), D_1(r_{D_1},\pi/2)$, where
\begin{eqnarray}
r_{C_1} & = & \exp\left(-\frac{\pi}{\sqrt{4m^2-1}}\right)\\
r_{D_1} & = & \exp\left[-\frac{1}{\sqrt{4m^2-1}}\left(\pi-\cot^{-1}\frac{1}{\sqrt{4m^2-1}}\right)\right]
\end{eqnarray}
These are the points with the largest radius along these axes, that can be reached from $A(1,0)$. Using (\ref{firstswitchingradius}), (\ref{kappa1}), (\ref{secswitchingradius}), (\ref{exponent}) and (\ref{secswitchingangle}) we find that the second switching curve $BC_1$ crosses the axis $\theta=\pi/2$ at the point $D_2(r_{D_2},\pi/2)$ where
\begin{equation}
r_{D_2}=\frac{\sqrt{m^2+2}}{1+\sqrt{m^2+1}}\exp\left[-\frac{1}{\sqrt{4m^2-1}}\cot^{-1}\left(\frac{m^2-1}{\sqrt{4m^2-1}}\right)\right]
\end{equation}

As we mentioned above, switching takes place only for $\kappa>\sqrt{m^2-1}$. For $\kappa\leq\sqrt{m^2-1}$ there is no switching and the optimal trajectory is given by (\ref{opttraj}). For these values of $\kappa$, the optimal trajectory crosses the axis $\theta=\pi,\pi/2$ at points with radius
\begin{equation}
r_{\pi}=\frac{\sqrt{\kappa^2+1}-1}{\sqrt{\kappa^2+1}+1},\,\,\,r_{\pi/2}=\frac{\kappa}{\sqrt{\kappa^2+1}+1}
\end{equation}
respectively, both increasing functions of $\kappa\geq 0$.
In Fig. \ref{fig:mgreater} we plot the optimal trajectory without switching for the largest permissible value $\kappa=\sqrt{m^2-1}$
\begin{equation}
\label{noswitching}
r_0(\theta_0)=\frac{\cos\theta_0+\sqrt{\cos^2\theta_0+m^2-1}}{m+1}
\end{equation}
and for $\theta_0\in[\theta_B,\pi]$. It crosses the axes $\theta=\pi,\pi/2$ at the points $C_2(r_{C_2},\pi), D_3(r_{D_3},\pi/2)$ where
\begin{equation}
r_{C_2}=\frac{m-1}{m+1},\,\,r_{D_3}=\frac{\sqrt{m^2-1}}{m+1}
\end{equation}
These are the largest radius points along these axes, that can be reached without switching.

Using the construction shown in Fig. \ref{fig:mgreater} we can find the switching points and the optimal control for any final point of the form $F_1(r_\tau,\pi)$ or $F_2(r_\tau,\pi/2)$. If $F_1\in C_1C_2$ then the optimal trajectory after the second switching is
\begin{equation}
\label{opttraj3}
r(\theta)=r_2(\theta_2)\frac{\cos\theta+\sqrt{\cos^2\theta+\kappa^2}}{\cos\theta_2+\sqrt{\cos^2\theta_2+\kappa^2}}
\end{equation}
where
\begin{equation}
\label{optkappa}
\kappa^2=(1-m\cot\theta_2)^2+m^2-1
\end{equation}
and $S_2(r_2,\theta_2)$ is the second switching point. Eq. (\ref{opttraj3}) is found after integrating (\ref{difeq}) from $S_2$ to $F_1$, while (\ref{optkappa}) holds because $S_2\in BC_1$, see (\ref{switchingangle}). The final point $F_1(r_\tau,\pi)$ belongs to this curve, so we find
\begin{equation}
\label{crossing2}
r_2(\theta_2)=r_\tau\frac{\cos\theta_2+\sqrt{\cos^2\theta_2+\kappa^2}}{-1+\sqrt{1+\kappa^2}}
\end{equation}
Plotting this curve for $\theta_2\in[\theta_B,\pi]$, it crosses the second switching curve at the second switching point $S_2$. The first switching angle can be found from
\begin{equation}
\label{firstswitching angle}
\theta_1=\cot^{-1}(2/m-\cot\theta_2).
\end{equation}
Having determined the two switching angles and the optimal $\kappa$, the optimal feedback control is also determined. If $F_1\in C_2O$ then no switching is necessary and the optimal control is given by (\ref{u}) with $\kappa=\kappa_{\pi}$ given by (\ref{kappa}).

For $F_2\in D_1D_2$ there is only one switching. The optimal trajectory after the switching is
\begin{equation}
\label{opttraj2}
r(\theta)=r_1(\theta_1)\sqrt{\frac{2m-\sin{2\theta_1}}{2m-\sin{2\theta}}}\exp{\left[-\frac{f(\theta_1,\theta)}{\sqrt{4m^2-1}}\right]}
\end{equation}
The final point $F_2(r_\tau,\pi/2)$ belongs to this curve, so we find
\begin{equation}
\label{crossing1}
r_1(\theta_1)=r_\tau\sqrt{\frac{2m}{2m-\sin{2\theta_1}}}\exp{\left[\frac{f(\theta_1,\pi/2)}{\sqrt{4m^2-1}}\right]}
\end{equation}
Plotting this curve for $\theta_1\in[0,\theta_B]$, it crosses the first switching curve at the switching point $S_1(r_1,\theta_1)$, so we find the switching angle and the corresponding optimal $\kappa$.
If $F_2\in D_2D_3$ then there are two switchings and the situation is similar to the one previously described, with
\begin{equation}
\label{crossing3}
r_2(\theta_2)=r_\tau\frac{\cos\theta_2+\sqrt{\cos^2\theta_2+\kappa^2}}{\kappa}
\end{equation}
instead of (\ref{crossing2}), since $\theta=\pi/2$ at the final point.
If $F_2\in D_3O$ then no switching is necessary; the optimal control is given by (\ref{u}) with $\kappa=\kappa_{\pi/2}$ given by (\ref{kappa}).

\subsubsection{$1/2<m<1$}

%% =================== Figure 5 ==================
\begin{figure}[t]
\centering
\includegraphics[scale=0.45]{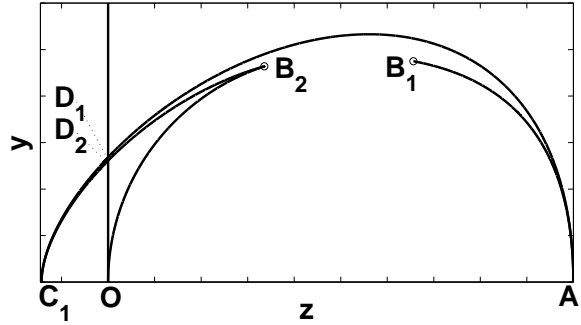}
\caption{Switching curves $AB_1$ and $B_2C_1$ for the case $1/2<m<1$. The outermost curve $AC_1$ defines the reachable set from the starting point $A(1,0)$ when $u\leq m$. Curve $B_2O$ is the optimal trajectory, after the second switching, for the lowest value $\kappa=0$. The perpendicular axis $\theta=\pi/2$ crosses the curves $AC_1, B_2C_1$ at the points $D_1, D_2$, respectively. Note that for these values of $m$, switching is inevitable for the transfers that we examine.}
\label{fig:mlower}
\end{figure}

For this case the two switching curves are described by the same equations as before but now there is no common point but a gap between them, due to the fact that the two roots of (\ref{switch}) are distinct for all $\kappa\geq 0$, see Fig \ref{fig:mlower}. The range for the switching angles is $\theta_1\in[0,\theta_{B_1}], \theta_2\in[\theta_{B_2},\pi]$ where
\begin{eqnarray}
\theta_{B_1} & = & \cot^{-1}\left(\frac{1+\sqrt{1-m^2}}{m}\right)\\
\theta_{B_2} & = & \cot^{-1}\left(\frac{1-\sqrt{1-m^2}}{m}\right)
\end{eqnarray}
The switching points $B_1,B_2$ correspond to the case $\kappa=0$ and their radius is
\begin{eqnarray}
r_{B_1} & = & \sqrt{\frac{1+\sqrt{1-m^2}}{2}}\\
r_{B_2} & = & r_{B_1}\times\exp{\left[-\frac{1}{\sqrt{4m^2-1}}\cot^{-1}\left(\frac{2m^2-1}{\sqrt{4m^2-1}\sqrt{1-m^2}}\right)\right]}
\end{eqnarray}
For the lowest value $\kappa=0$ we plot the optimal curve after the second switching. It is
\begin{equation}
\label{0switching}
r_0(\theta_0)=r_{B_2}\frac{\cos\theta_0}{\cos\theta_{B_2}}
\end{equation}
This curve meets the axis $\theta=\pi/2$ at the origin $O$.
Points $C_1, D_1, D_2$ are defined as before and their coordinates are the same. The situation is as depicted in Fig. \ref{fig:mlower}. Observe that for $1/2<m<1$, switching is inevitable for the transfers that we examine; since the control is more restricted than the previous case, the boundary value has to be used to achieve the desired transfers. Any final point of the form $F_1(r_\tau,\pi)$ belongs to $C_1O$ and there are two switchings in the corresponding optimal trajectory. The same happens for $F_2(r_\tau,\pi/2)\in D_2O$, while for $F_2\in D_1D_2$ there is only one switching.

\subsubsection{Summary of the results }

We summarize the above results
\begin{itemize}
 \item $m\geq 1$
  \begin{enumerate}
   \item $\theta=\pi$
    \begin{enumerate}
     \item $r_{C_2}<r_\tau\leq r_{C_1}$, two switchings
     \item $r_\tau\leq r_{C_2}$, no switching
    \end{enumerate}
   \item $\theta=\pi/2$
    \begin{enumerate}
     \item $r_{D_2}<r_\tau\leq r_{D_1}$, one switching
     \item $r_{D_3}<r_\tau\leq r_{D_2}$, two switchings
     \item $r_\tau\leq r_{D_3}$, no switching
    \end{enumerate}
  \end{enumerate}
 \item $1/2<m<1$
  \begin{enumerate}
   \item $\theta=\pi$
    \begin{enumerate}
     \item $r_\tau\leq r_{C_1}$, two switchings
    \end{enumerate}
   \item $\theta=\pi/2$
    \begin{enumerate}
     \item $r_{D_2}<r_\tau\leq r_{D_1}$, one switching
     \item $r_\tau\leq r_{D_2}$, two switchings
    \end{enumerate}
  \end{enumerate}
\end{itemize}

If $\theta_1,\theta_2$ are the switching angles then the optimal control for $\theta\in(\theta_1,\theta_2)$ is
\begin{equation}
\label{optin}
u(\theta)=m
\end{equation}
while outside this interval is
\begin{equation}
\label{optout}
u(\theta)=\sin\theta(\cos\theta+\sqrt{\cos^2\theta+\kappa^2})
\end{equation}
where
\begin{equation}
\label{kappaopt}
\kappa=\sqrt{(m\cot\theta_1-1)^2+m^2-1}=\sqrt{(1-m\cot\theta_2)^2+m^2-1}.
\end{equation}
The switching angles are calculated as described above and are related through
\begin{equation}
\label{thitas}
\cot\theta_1+\cot\theta_2=2/m.
\end{equation}

\section{Examples}

\begin{figure}[t]
 \centering
		\begin{tabular}{cc}
     	    \subfigure[$\ $$m=2,F_1(0.39,\pi)$]{
	            \label{fig:a}
	            \includegraphics[width=.45\linewidth]{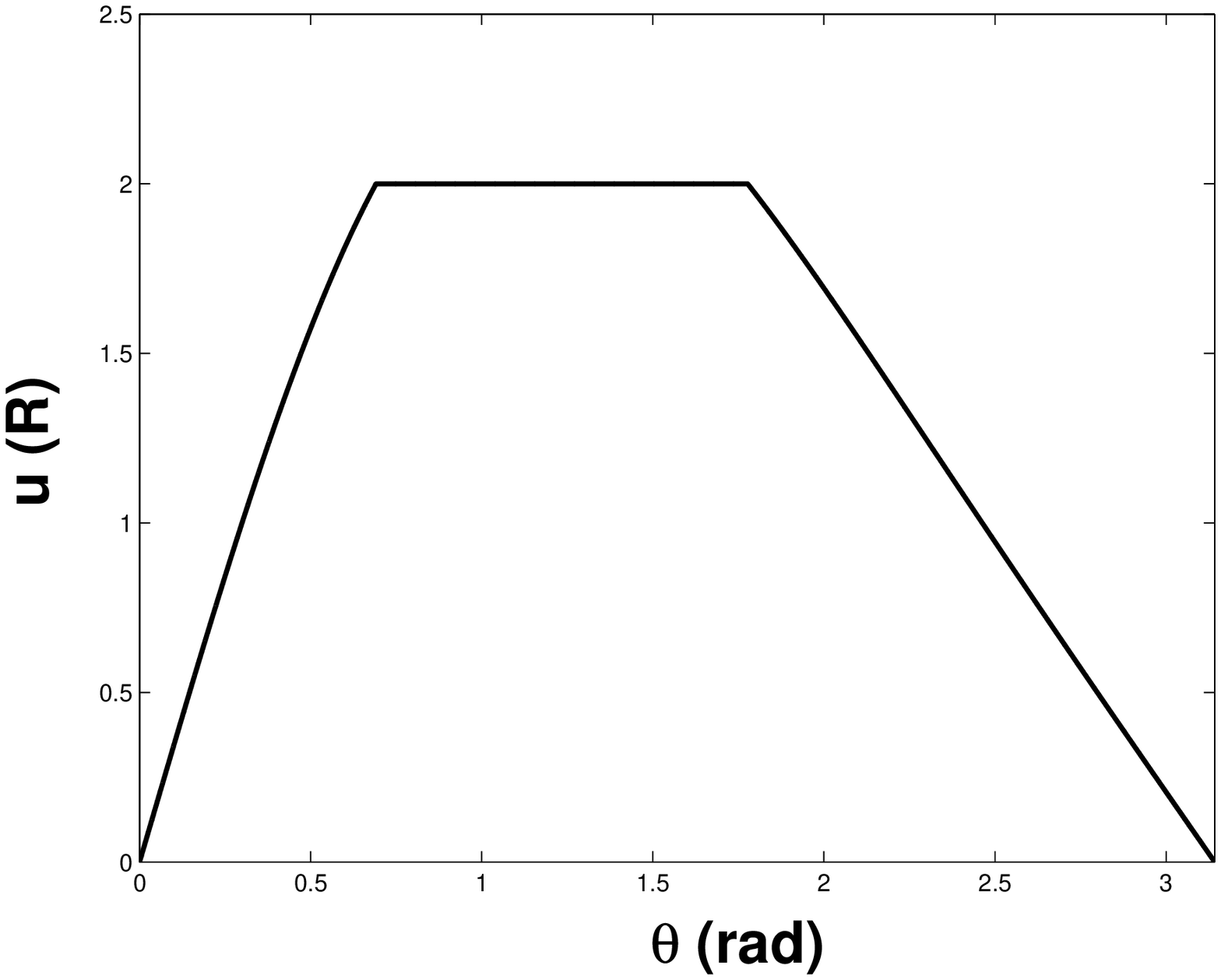}} &
	        \subfigure[$\ $Optimal trajectory (dashed line)]{
	            \label{fig:b}
	            \includegraphics[width=.45\linewidth]{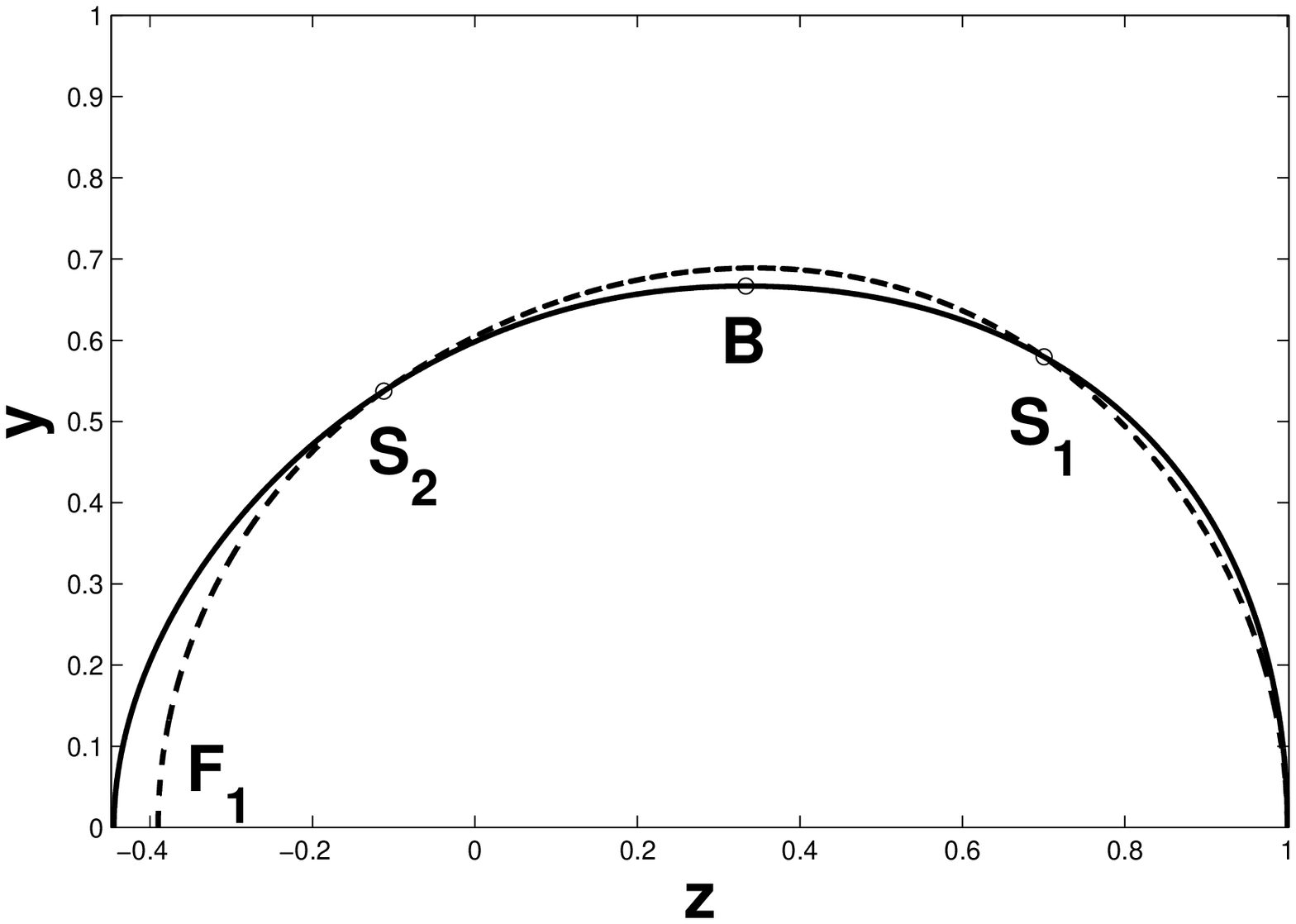}}\\
            \subfigure[$\ $$m=2,F_2(0.61,\pi/2)$]{
	            \label{fig:c}
	            \includegraphics[width=.45\linewidth]{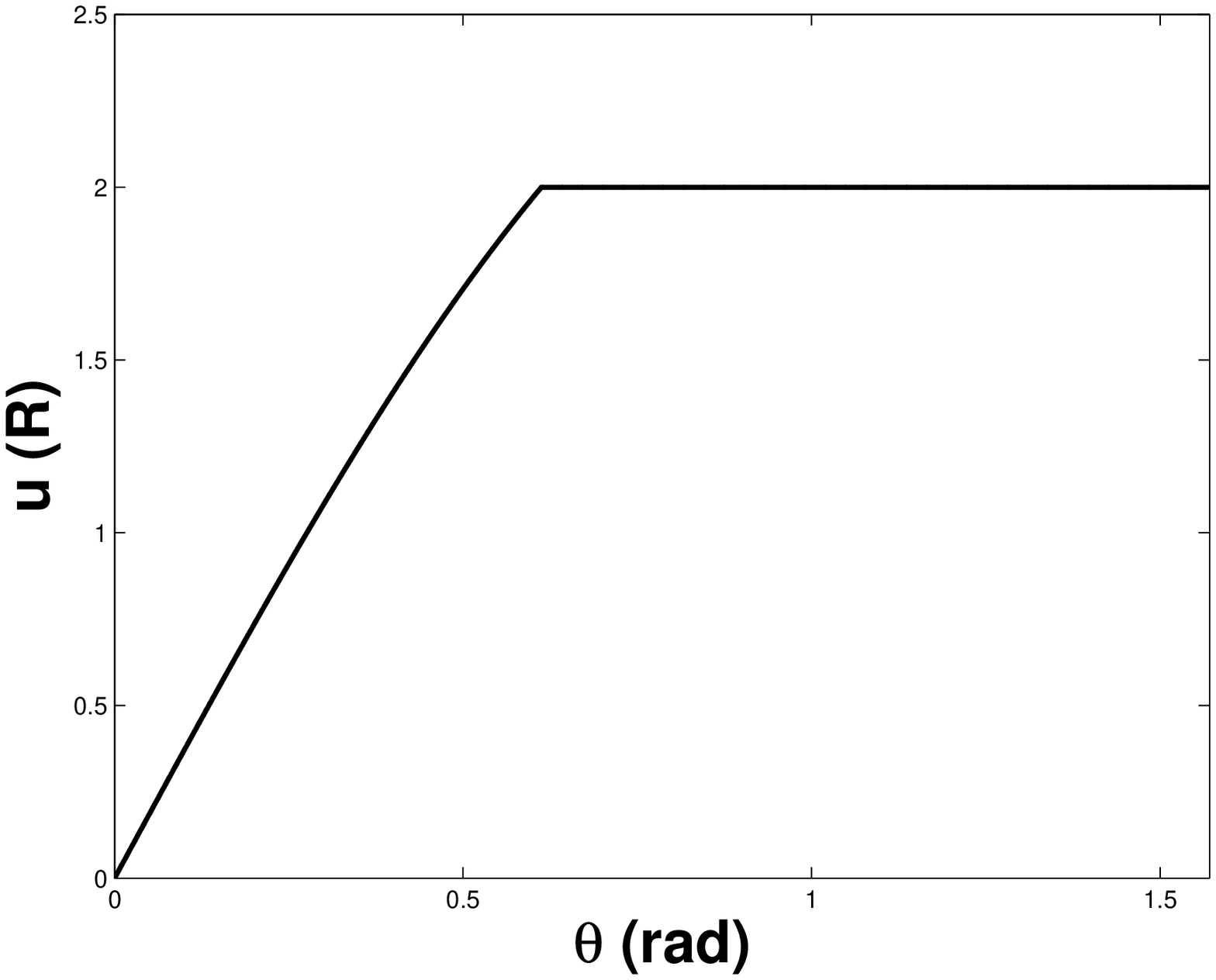}} &
	        \subfigure[$\ $Optimal trajectory (dashed line)]{
	            \label{fig:d}
	            \includegraphics[width=.45\linewidth]{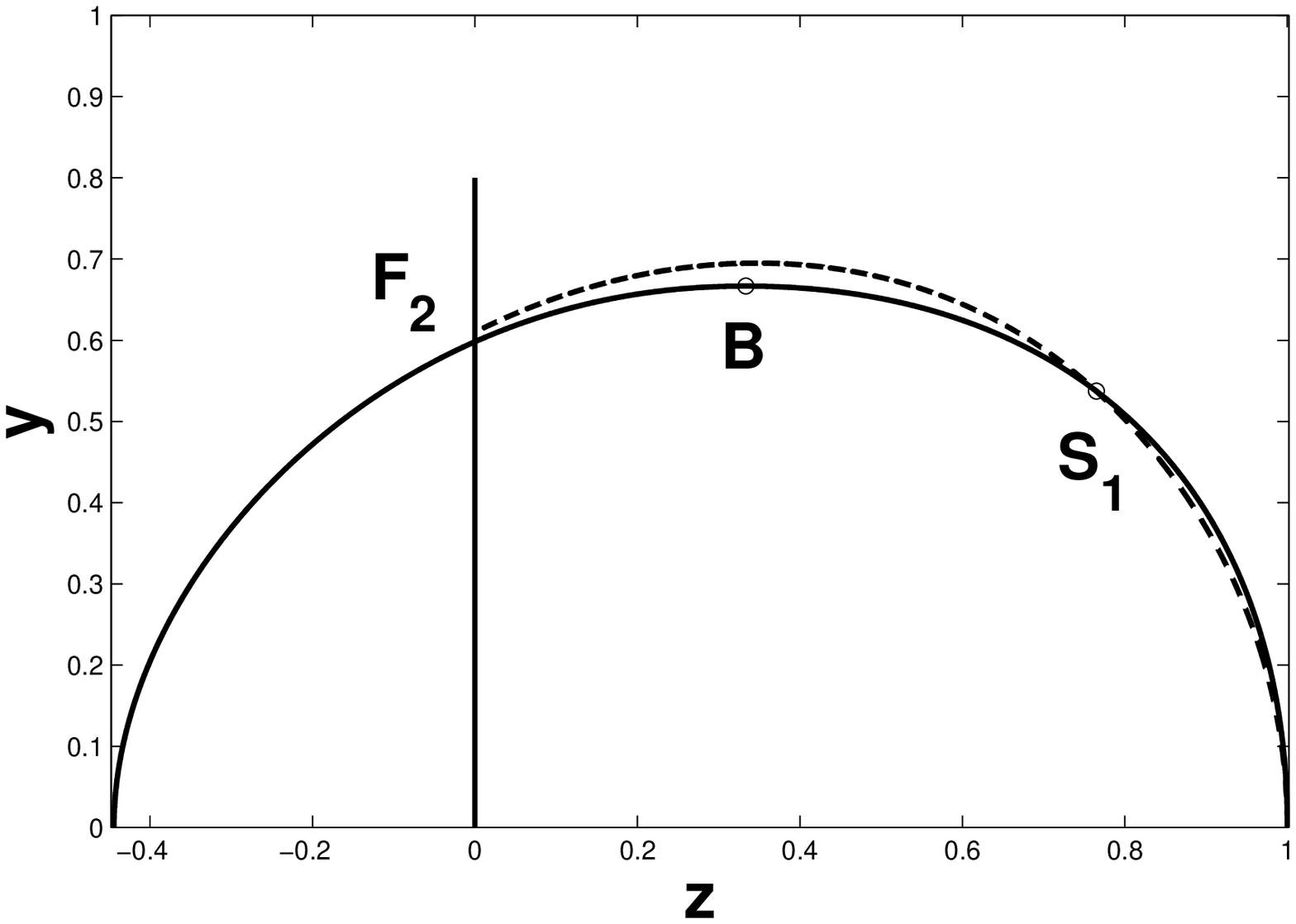}}\\
            \subfigure[$\ $$m=0.95,F_2(0.2,\pi/2)$]{
	            \label{fig:e}
	            \includegraphics[width=.45\linewidth]{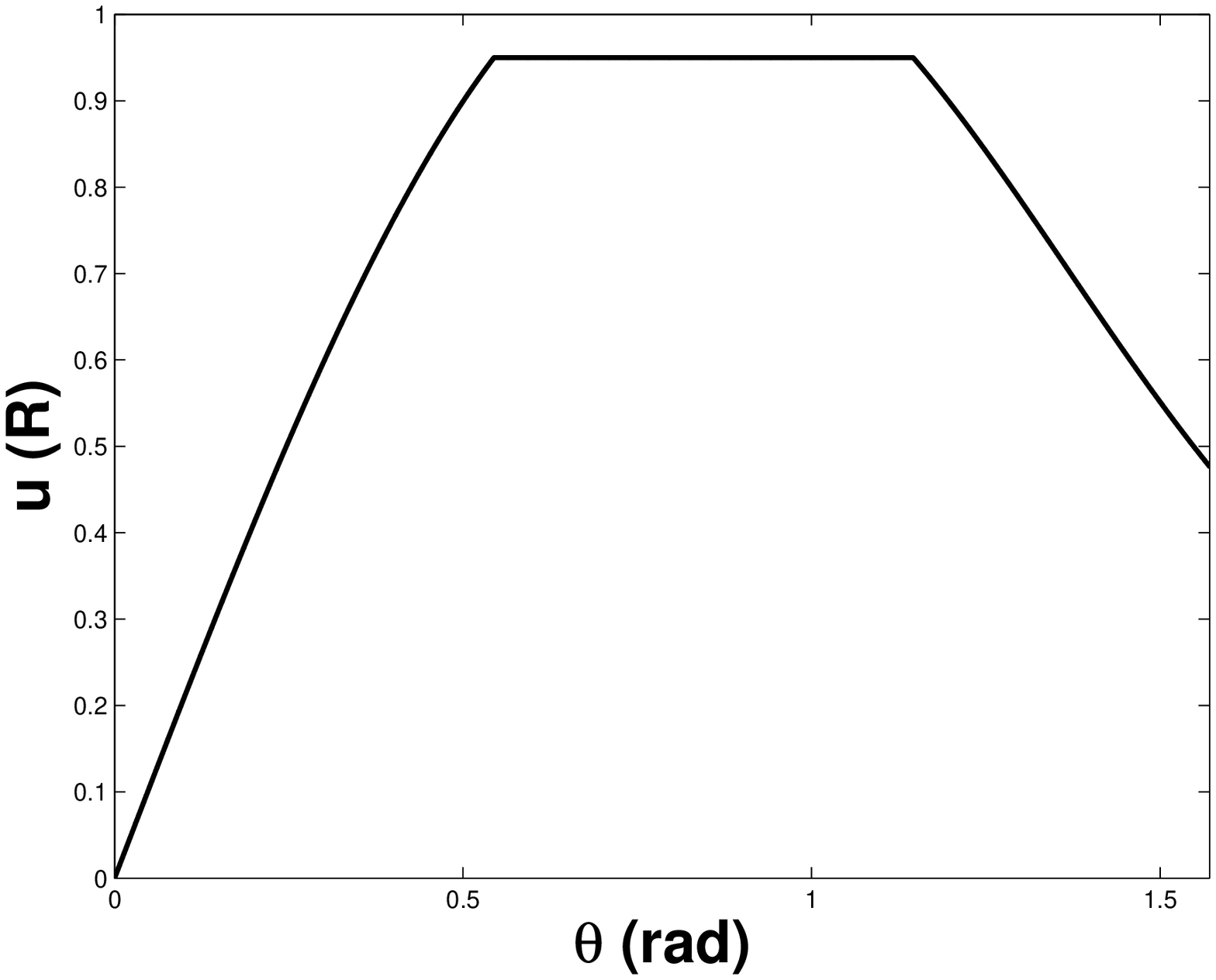}} &
	        \subfigure[$\ $Optimal trajectory (dashed line)]{
	            \label{fig:f}
	            \includegraphics[width=.45\linewidth]{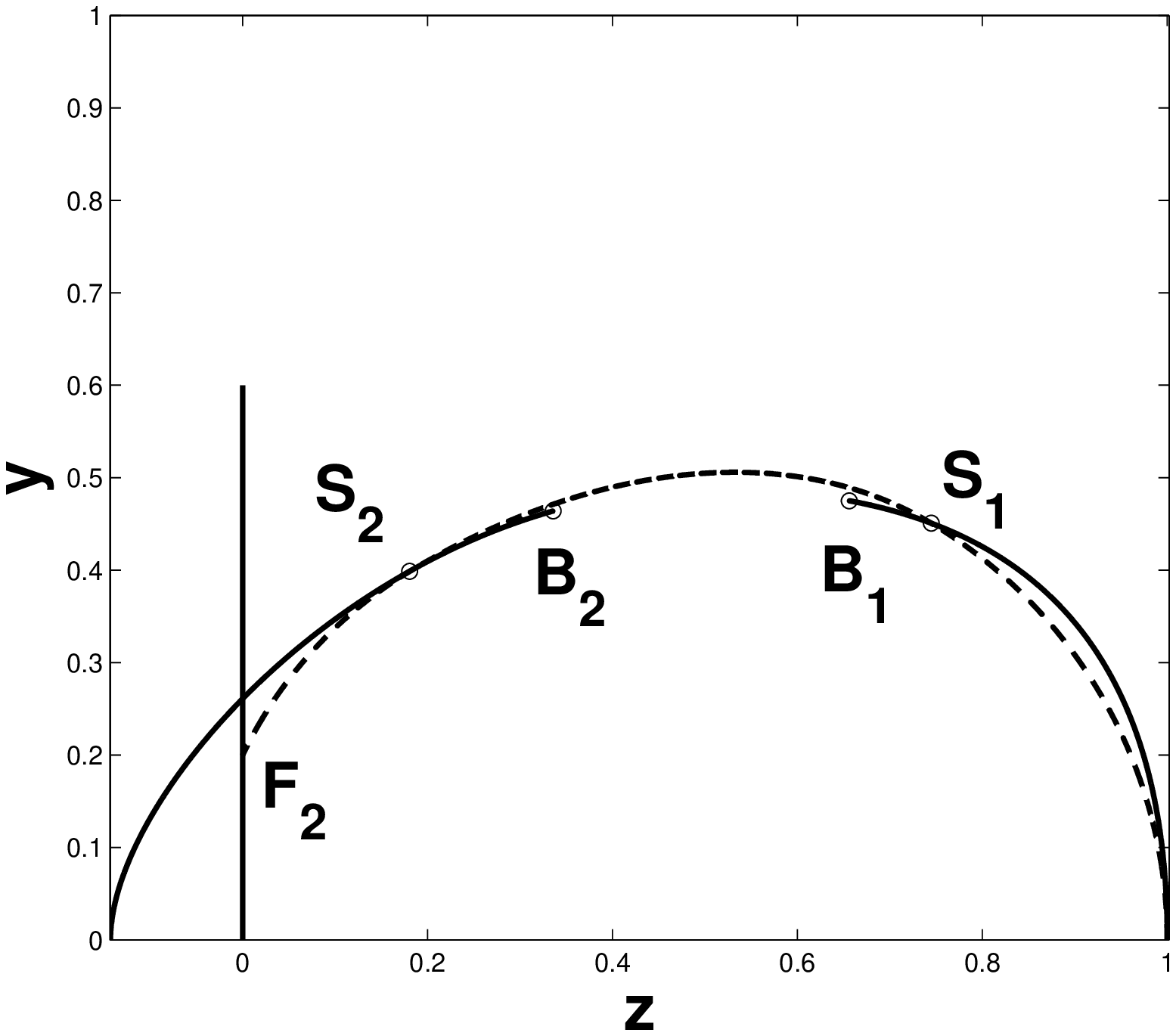}}
		\end{tabular}
\caption{minimum-energy pulses (a, c, e) and corresponding optimal trajectories (b, d, f) for various values of the control upper bound $m$ and the coordinates $(r_\tau,\theta_\tau)$ of the final point. The switching curves and points are also shown.}
 \label{fig:pulses1}
\end{figure}

Here we  present some examples using specific values for the upper bound $m$ and the coordinates $(r_\tau,\theta_\tau)$ of the final point. We start with the case $m=2,r_\tau=0.39,\theta_\tau=\pi$. Using the results of the previous section, it is not hard to see that the optimal trajectory contains two switching points. Plotting Eq. (\ref{crossing2}) we find that it crosses the second switching curve at the switching point $S_2$ with angle $\theta_2=1.7766$ rad. The first switching angle is found from (\ref{thitas}) to be $\theta_1=0.6912$ rad. The optimal control is given by (\ref{optin}), (\ref{optout}) with $\kappa=2.2382$, as it is determined from (\ref{kappaopt}). In Fig. \ref{fig:a}, \ref{fig:b} we plot the optimal pulse and the corresponding optimal trajectory. The two switching points $S_1,S_2$ are also shown.

For the next example we use the values $m=2,r_\tau=0.61,\theta_\tau=\pi/2$. In this case, the optimal trajectory contains only one switching point. Plotting Eq. (\ref{crossing1}) we find that it crosses the first switching curve at the switching point $S_1$ with angle $\theta_1=0.6124$ rad. The optimal control is given by (\ref{optout}) with $\kappa=2.5322$ for $\theta\in [0,\theta_1]$ and by (\ref{optin}) for $\theta\in (\theta_1,\pi/2]$. In Fig. \ref{fig:c}, \ref{fig:d} we plot the optimal control and the corresponding optimal trajectory, with the one switching point $S_1$.

The last case that we consider is $m=0.95,r_\tau=0.2,\theta_\tau=\pi/2$. The optimal trajectory contains two switching points. Plotting Eq. (\ref{crossing3}) we find that it crosses the second switching curve at the switching point $S_2$ with angle $\theta_2=1.1456$ rad. The first switching angle is found from (\ref{thitas}) to be $\theta_1=0.5442$ rad. The optimal control is given by (\ref{optin}), (\ref{optout}) with $\kappa=0.4766$, as it is determined from (\ref{kappaopt}). In Fig. \ref{fig:e}, \ref{fig:f} we plot the optimal pulse and the corresponding optimal trajectory. The two switching points $S_1,S_2$ are also shown.

%% =================== Figure 6 ==================

\section{Conclusion}

To conclude, in this letter we calculated minimum-energy $\pi/2$ and $\pi$ pulses for Bloch equations in the case where transverse relaxation dominates and the control amplitude is bounded, using optimal control theory. This work is expected to find applications in NMR Spectroscopy, Magnetic Resonance Imaging (MRI) and Quantum Information Processing, serving as a reference for numerical studies of more complicated and realistic situations that incorporate for example longitudinal relaxation and magnetic field inhomogeneity.

%\section{Acknowledgements}


\end{document}